\documentclass[fleqn,reqno,11pt]{amsart}
\usepackage[pdftex,marginratio=1:1,width=130mm,height=190mm,top=40mm,headsep=17mm,head=7mm,footskip=11.29mm]{geometry}
\usepackage[colorlinks,linkcolor=blue,citecolor=magenta,anchorcolor=blue,bookmarksopen,pdfpagetransition={Wipe}]{hyperref}
\usepackage{graphicx}
\usepackage{lineno,hyperref}
\usepackage{amsmath, amsthm, amscd, amsfonts, amssymb,graphicx, color, textcomp, mathtools}
\usepackage{pgf,tikz,pgfplots}
\pgfplotsset{compat=1.15}
\usepackage{mathrsfs}
\usetikzlibrary{arrows}
\usepackage{multirow}
\usepackage{xspace}
\usepackage{array}
\usepackage{multirow}
\usepackage{float}

\newcommand{\dx}{\,\mathrm{d{\bf x}}}
\DeclareSymbolFont{rmlargesymbols}{OMX}{mdbch}{m}{n}
\DeclareMathSymbol{\rmintop}{\mathop}{rmlargesymbols}{82}
\newcommand*{\Scale}[2][4]{\scalebox{#1}{$#2$}}%
\newcommand{\intt}{\Scale[1]{\rmintop\nolimits\kern-1pt}}
\newcommand{\B}{\bf x}
\newcommand{\bu}{\bf u}
\newcommand{\bt}{\boldsymbol{\theta}}
\modulolinenumbers[5]

\theoremstyle{definition}

\numberwithin{equation}{section}

\begin{document}
\unitlength 1cm\begin{picture}(0,0)\put(2,1){ \scriptsize
\textbf{(CMDE)} \textsl{Computational Methods for Differential
Equations}}\end{picture}

\title[Semi-Algebraic Mode Analysis For Finite Element Discretisations]{Semi-Algebraic Mode Analysis For Finite Element Discretisations Of The Heat Equation}%
\author[N. Habibi]{Noora Habibi$^*$}%
\address{Faculty of Applied Mathematics, Shahrood University Of Technology, P.O. Box 3619995161 Shahrood(Iran).}%
\email{habibi85nu@gmail.com}%
\thanks{$*$ corresponding}

\author[A. Mesforoush]{Ali Mesforoush}%
\address{Faculty of Applied Mathematics, Shahrood University Of Technology, P.O. Box 3619995161 Shahrood(Iran).}%
\email{ali.mesforush@shahroodut.ac.ir}%

\author[F. Gaspar J. L.]{Francisco J. Gaspar Lorenz}%
\address{Department of Applied Mathematics, Science Faculty, University of Zaragoza, Pedro Cerbuna, 12, 50009 Zaragoza(Spain)}%
\email{fjgaspar@unizar.es}%

\author[C.  Rodrigo]{ Carmen Rodrigo}%
\address{Department of Applied Mathematics,  School of Engineering and Architecture, University of Zaragoza, Maria de Luna, 3, 50018, Zaragoza(Spain)}%
\email{carmenr@unizar.es}%


\begin{abstract}
In this work, a semialgebraic mode analysis (SAMA) is proposed for  investigating the convergence of a multigrid waveform relaxation method applied to the Finite Element (FE) discretization of the heat equation in two and three dimensions. This analysis for finite element methods is more involved and more general  than that for Finite Difference (FD) discretizations, since mass matrix must be considered. The proposed analysis results in a very useful tool to study the behaviour of the multigrid waveform relaxation method depending on the parameters of the problem. 
\end{abstract}
\maketitle
\section{Introduction}\label{sec1}
The Semi-Algebraic Mode Analysis (SAMA) was introduced by Friedhoff and MacLachlan in \cite{1} for predicting the convergence factor of time-dependent Partial Differential Equations (PDEs) as a generalization of Local Fourier Analysis (LFA).  LFA is the most powerful technique for predicting the convergence of multigrids methods \cite{2, 3}. This method analyzes the behavior of the local components involved in multigrid methods on a basis of complex exponential functions.  LFA, however, has not been successful in predicting the convergence factor for time dependent PDE problems   \cite{7, 1},  and to overcome this failure, SAMA was introduced.   The main idea of SAMA is the combination of LFA only in space with an exact analytical approach in time. In this work, we utilize SAMA to analyze the convergence factor of multigrid waveform relaxation applied to finite element discretization of the heat equation. 


Let us consider as model problem the heat equation with homogeneous Dirichlet boundary conditions
\begin{equation}\label{1.1}
\begin{split}
&D_tu({\B},t)-\Delta u({\B},t))= f({\B},t),\quad {\B} \in \Omega,\quad t>0,\\
&u({\B},t)=0,\quad \text{on } \partial \Omega,\quad t>0,\\
&u({\B},0)=g({\B}), \quad {\B}\in \Omega,
\end{split}
\end{equation}
where $\Omega \subset \mathbb{R}^d$,  for $d=2$ (or $d=3$) is a bounded domain with boundary $\partial \Omega$. 
In order to establish the finite element approximation of our problem, let $\Omega_h$ be a triangulation of $\Omega$, satisfying the usual admissibility
assumption, i.e. the intersections of two different elements is either empty, a vertex, or a whole edge.
%
%
%
%
  Let $V_h$ be the finite element space of continuous piecewise bilinear (or trilinear) functions associated with $\Omega_h$ vanishing on the boundary $\partial \Omega$. The discrete approximation $u_h\in V_h$  solves the following problem
 \begin{align*}
\left(D_tu_h,v_h\right)+a\left(u_h,v_h\right)=(f,v_h)\quad  v_h\in V_h,
\end{align*}
where
\begin{align*}
&\left(D_tu_h,v_h\right) = \int_{\Omega} (D_tu_h) v_h\, \dx,\quad  \quad (f,v_h) = \int_{\Omega} fv_h \, \dx\\
&a\left(u_h,v_h\right) = \int_{\Omega}\nabla u_h \cdot \nabla v_h\, \dx.
\end{align*}
Let \{$\phi_1, \ldots ,\phi_{N}$\} be the nodal basis of $V_h,$ i.e., $\phi_i({\B}_j) = \delta_{ij}$, with ${\B}_j$ an interior node of the mesh ${\Omega}_h$. 
The approximation $u_h=\sum_{i=1}^{N} u_i(t)\phi_i({\B})$ is found by solving the following set of equations,
\begin{align*}
\left(D_tu_h,\phi_j\right)+a\left(u_h,\phi_j\right)=(f,\phi_j),\quad \text{for }j=1,2,\ldots,N.
\end{align*}
We rewrite these equations in terms of the mass matrix $B=\{(\phi_i,\phi_j)\}$ and the stiffness matrix $A=\{a(\phi_i,\phi_j)\}$, in a more standard form, as a system of ordinary differential equations (ODEs)
\begin{align}\label{1.2}
B_h\dot{{\bu}}_h(t)+A_h{\bu}_h(t)=F_h(t),\quad {\bu}_h(0)=g_h,\quad t>0,
\end{align}
where ${\bu}_h(t)=[u_1(t),u_2(t),\ldots,u_{N}(t)]^t\in \mathbb{R}^{N}$ and the coefficient matrices $A_h, B_h\in \mathbb{R}^{N\times N}$ and the right hand side  $F_h(t)=[(f,\phi_1),(f,\phi_2),\ldots,$ $(f,\phi_{N})]^T$ $ \in \mathbb{R}^{N}$ are considered.
\begin{figure}[!ht]
\centering
	\centerline{\includegraphics[width=11cm]{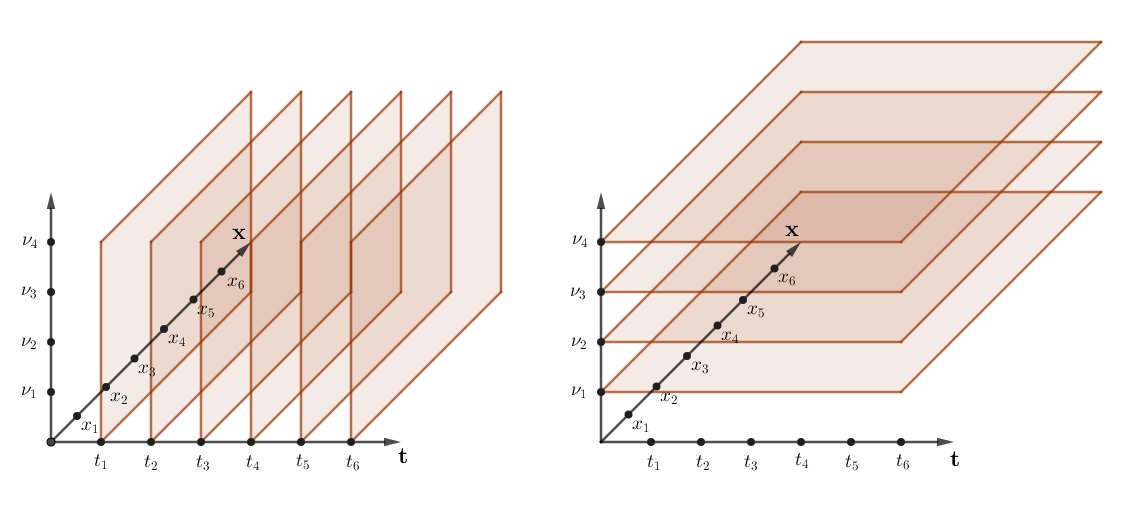}}
\caption{Time-stepping (at left) versus waveform relaxation (at right). In this figure  $\nu_i,$ for $i=1,2,\ldots $ indicates to the iterations of two methods.}
\label{TW1}
\end{figure}

We have many choices to pick a suitable method for solving the obtained ODE system \eqref{1.2}. In general, we can divide all the well known methods into two classes: time-marching approaches and time parallel techniques. 
In a time marching approach, we  solve in each time step a time-independent problem then we go to the next time step. Indeed, as we can get out from its name, each time step will be solved after the other in a sequential manner,
 see Figure \ref{TW1} (left). Although this approach is simpler,  when we need multiprocessing capability or parallelization of the temporal variable we have to seek for another approach. So we consider a full space-time method.
  The time parallel class itself can be classified into four groups of methods (see \cite{8}): Multiple shooting; Domain decomposition and waveform relaxation; Space-time multigrid method; and Direct time parallel methods. Here, we only concentrate on the waveform relaxation method and in  particular, on its multigrid extension.
  
The waveform relaxation method (WR) is a technique for solving ordinary differential equations. It can be also applied to time-dependent PDEs when their spatial derivatives are replaced by a discrete formula (obtained by the Finite Element method in our particular  case) as viewpoint of the method of lines scheme. The WR method is based on splitting matrices $A_h$ and $B_h$ as $B_h=M_{B_h}-N_{B_h}$ and $A_h=M_{A_h}-N_{A_h}$, leading to the following iteration
\begin{align}\label{eq1}
M_{B_h}\dot{{\bu}}_h^{k}(t)+M_{A_h}{\bu}_h^{k}(t)=N_{B_h}\dot{{\bu}}_h^{k-1}(t)+N_{A_h}{\bu}_h^{k-1}(t)+F_h(t),
\end{align}
where $ {\bu}_h^k(0)=g_h,$ for $k\geq 1$ and ${\bu}_h^k(t)$ indicates the approximation of ${\bu}(t)$ at iteration $k$. It is natural to define ${\bu}_h^0(t)$ along the whole time interval equal to the initial condition, i.e. ${\bu}_h^0(t)=g_h,\,\, t>0$. 
Considering the decomposition of matrices $A_h$ and $B_h$  as $A_h=-L_{A_h}+D_{A_h}-U_{A_h}$ and $B_h=-L_{B_h}+D_{B_h}-U_{B_h}$, where $L_{A_h}$ and $L_{B_h}$ are strictly lower triangular matrices, $D_{A_h}$ and $D_{B_h}$ are diagonal matrices, and $U_{A_h}$ and $U_{B_h}$ are strictly upper triangular matrices,  for the Gauss-Seidel waveform relaxation method, which is considered in this work, the splittings in \eqref{eq1} are as follows:
\begin{align*}
\begin{matrix}
M_{A_h}=-L_{A_h}+D_{A_h}  &  N_{A_h}=U_{A_h},\\
M_{B_h}=-L_{B_h}+D_{B_h}  &  N_{B_h}=U_{B_h}.
\end{matrix}
\end{align*}

We use multigrid technique to accelerate the convergence of the Gauss-Seidel waveform relaxation method. A multigrid acceleration of this method was firstly studied by Lubich and Ostermann in \cite{9} and independently developed in \cite{10}. 

 In order to apply a geometric multigrid waveform relaxation procedure the coarsening applies only in the spatial domain  and  we consider a hierarchy of grids like $\Omega_{2^lh}\subset \ldots \subset$
$\Omega_{2h}\subset\Omega_{h}$.  We obtain a new iterate ${\bu}_h^{(k)}$ from the former waveform ${\bu}_h^{(k-1)}$ in three steps: Pre-smoothing, coarse grid correction and post smoothing. In Algorithm 1 we present the multigrid waveform relaxation algorithm (WRMG) depending on the defined Gauss-Seidel  waveform relaxation method as smoother and the rest of the operators involved in the multigrid procedure. We consider  standard coarsening for constructing the coarse meshes and discretization coarse grid approximation (DCA) in coarser grids. Regarding intergrid transfer operators, the interpolation operators are the nine point stencil operator corresponding to the bilinear interpolation for two dimensional problems and its generalization for three dimensional problems, that is,  the  trilinear interpolation operator. The restriction operators are considered as the adjoint of the prolongation operators.
\begin{table}[ht]
\small
\begin{tabular}{l}
  \hline
  {\bf Algorithm 1}  Multigrid waveform relaxation based on Gauss-Seidel smoother\\
  $ {\bu}_h^{(k)}(t)\rightarrow  {\bu}_h^{(k+1)}(t)$ \\
  \hline
  {\bf If} we are on the coarsest grid (given by spatial grid size $2^lh=h_0$), then solve\\
   the following   equation by a direct or fast solver\\
  $\hspace{1cm}B_{h_0}\dot{{\bu}}_{h_0}^{k+1}(t)+A_{h_0}{\bu}_{h_0}^{k+1}(t)=F_{h_0}(t)$.\\
  {\bf Else}\\
  (\textit{Presmoothing})  Perform $k_1$ steps of Gauss-Seidel waveform relaxation, \\
  $\hspace{1cm}v_h^k(t)=S^{k_1}\left({\bu}_h^{k}(t) \right)$.\\
  (\textit{Coarse grid correction})\\
 $\hspace{.88cm}$Compute the defect\\
   $\hspace{1cm}\bar{d}_h^k(t)=F_h(t)-B_h\dot{v}_h^{k}(t) - A_hv^{k}_h(t)$\\
  $\hspace{1cm}$Restrict the defect\\ 
  $\hspace{1cm}\bar{d}_{2h}^k(t)=I_h^{2h}\bar{d}_h^k(t)$\\
  $\hspace{1cm}$Perform $\gamma\geq 1$ cycles of WRMG on $\Omega_{2h}$ to solve the following defect \\
  $\hspace{1cm}$equation,\\
   $\hspace{1cm}B_{2h}\dot{{e}}_{2h}^k(t)+A_{2h}{e}_{2h}^k(t)=\bar{d}_{2h}^k(t),\,\,{e}_{2h}^k(0)=0$\\
   $\hspace{1cm}$Interpolate the correction\\
   $\hspace{1cm}{e}_h^k(t)=I_{2h}^h{e}_{2h}^k(t)$\\
  $\hspace{.9cm}$Correct the current approximation with the interpolation of the correction, \\
  $\hspace{1cm}v_h^{k+1}(t)=v_h^{k}(t)+{e}_{h}^k(t)$. \\
   (\textit{Postsmoothing}) Perform $k_2$ steps of Gauss-Seidel waveform relaxation, \\
  $\hspace{1cm}{\bu}^{k+1}_h(t)=S^{k_2}\left(v_h^{k+1}(t) \right)$.\\
  {\bf End If}\\
  \hline
\end{tabular}
\end{table} 
 
In Algorithm 1,  using the Crank-Nicolson approach for time discretization we obtain a space-time multigrid method with coarsening only in space. Thus, we have time-line Gauss-Seidel waveform relaxation, with standard full weighting restriction and bilinear interpolation in space for data transfer between the levels in the multigrid hierarchy.

In the rest of this work,  first we explain the two-grid SAMA analysis in two dimensions for Gauss-Seidel waveform relaxation with multigrid acceleration in Section \ref{s2.1}. At the end of this section , Section \ref{sec2.2}, we will present our analysis for the considered model problem \eqref{1.2}, by defining all required stencils involved in the FE method with bilinear basis functions. Section 3 is devoted to extending  SAMA analysis to three dimensions. Again, at the end of Section 3 we perform some analysis for our particular model problem introducing all necessary stencils for the FE method with trilinear basis functions. Conclusions are drawn in Section 4.
\section{SAMA in two dimensions}
%
\subsection{Theoretical framework}\label{s2.1}
Now, we describe the convergence analysis of the multigrid waveform relaxation method by SAMA.  For this end, we will  do  the following steps. 
\begin{enumerate}
\item 
First, we present some primary definitions,
\item Next, we explain SAMA smoothing analysis for the Gauss-Seidel relaxation procedure,
\item Then, we investigate the analysis for the coarse-grid correction operator,
\item Finally, we combine the two last steps to perform a complete two-grid analysis.
\end{enumerate}

 We define the infinite grid
$\mathcal{Q}_h=\{ {\bf x=kh}=(k_1h_1,k_2h_2), \,\, {\bf k}\in\mathbb{Z}^2\},$
and the so-called Fourier modes as $\varphi_h(\bt, {\B})=e^{i\bt\cdot {\B}}=e^{i\theta_1x_1}e^{\theta_2x_2}$ where ${\bf h}$ is the spatial discretization step and        $\bt \in \Theta_h=(-\pi/h,\pi/h]^2$. Now, we can define any discrete grid function for a fixed t as  a formal linear combination of the Fourier modes by 
\begin{align}\label{eq4.1}
u_h({\B},t)=\sum_{\bt\in \Theta_h}c_{\bt}(t)\varphi_h(\bt,{\B}),\quad {\B}\in \mathcal{Q}_h,
\end{align}
where coefficients $c_{\bt}(t)$  depend on the time variable. The  Fourier modes yield the so-called Fourier space $\mathcal{F}(\mathcal{Q}_h)=\{\varphi_h(\bt,{\B}), \,\bt \in \Theta_h\}$, 
and they are formal eigenfunctions of any discrete operator $L_h$, e.g. for  the standard discrete Laplace operator,   $L_h=-\Delta_h=\frac{1}{h^2}[-1\,\,2\,\,-1]$,  the expression $L_h\varphi_h(\bt,{\B})=\widehat{L}_h({\bt})\varphi_h(\bt,{\B})$
 holds where
\[\widehat{L}_h({\bt})=\frac{2}{h^2}\left(1-\cos(\theta h)\right),\]
 is the Fourier representation of $L_h$ on the Fourier space, also called formal eigenvalue or the Fourier symbol of $L_h$.

Now,   we study the effect of Gauss-Seidel relaxation procedure and the coarse-grid correction operator on the Fourier modes.
To reach this goal, we need to define the high- and low-frequency components with respect to the coarsening strategy. We consider $\Theta_{2h}=(-\pi/2h,\pi/2h]^2$ and $\Theta_h\smallsetminus\Theta_{2h}$ as the low- and high frequency spaces, respectively. Here, we have used the standard coarsening which means that the step size is double on the coarse grid, denoted by $\mathcal{Q}_{2h}$.

By considering $A_h=M_{A_h}-N_{A_h}$ and $B_h=M_{B_h}-N_{B_h}$ that have been defined in the previous section, we can write an iteration of the waveform relaxation method for an error grid function as:
\begin{equation}
M_{B_h}\dot{e}_h^k({\B},t)+M_{A_h}e^k_h({\B},t)=N_{B_h}\dot{e}^{k-1}_h({\B},t)+N_{A_h}e^{k-1}_h({\B},t),  
\end{equation}
 for $k\geq1,\, {\B}\in\mathcal{Q}_h$ and $t>0$, where, $e^{k-1}_h(\cdot,t)$ and $e^{k}_h(\cdot,t)$ are the error grid functions at the $k-1$ and $k$ iterations  and $e_h^k({\B},0)=0$.

Following  equation \eqref{eq4.1}, we can write $e_h^i({\B},t)$ in the $i$th iteration as:
\begin{align}
e_h^i({\B},t)=\sum_{\bt\in \Theta_h}c_{\bt}^i(t)\varphi_h(\bt,{\B}),\quad {\B}\in \mathcal{Q}_h,\,\, t>0.
\end{align}
We denote as $\widehat{M}_{A_h}({\bt})$, $\widehat{M}_{B_h}({\bt})$, $\widehat{N}_{A_h}({\bt})$ and $\widehat{N}_{B_h}({\bt})$ the symbols of $M_{A_h}$, $M_{B_h}$, $N_{A_h}$ and $N_{B_h}$, respectively. So, for each frequency $\bt \in \Theta_h$ we have 
 \begin{align}\label{eq4.4}
 \widehat{M}_{B_h}(\bt)\dot{c}_{\bt}^k(t)+\widehat{M}_{A_h}(\bt)c_{\bt}^k(t)=\widehat{N}_{B_h}(\bt)\dot{c}_{\bt}^{k-1}(t)+\widehat{N}_{A_h}(\bt)c_{\bt}^{k-1}(t),
 \end{align}
for $ k\geq 1,\,\, t>0$. Considering a uniform grid in time with M subdivisions and time step size $\tau$, we apply the Crank-Nicolson time discretization on equation \eqref{eq4.4}, given by
\begin{equation}
\begin{split}
&\hspace{-1cm}\widehat{M}_{B_h}(\bt)\dfrac{c_{\bt,i}^k-c_{\bt,i-1}^k}{\tau}-\widehat{N}_{B_h}(\bt)\dfrac{c_{\bt,i}^{k-1}-c_{\bt,i-1}^{k-1}}{\tau}\\
&\hspace{-1cm}=\frac{1}{2}\left( -\widehat{M}_{A_h}(\bt)c_{\bt,i}^{k}+\widehat{N}_{A_h}(\bt)c_{\bt,i}^{k-1} \right)+
  \frac{1}{2}\left( -\widehat{M}_{A_h}(\bt)c_{\bt,i-1}^{k}+\widehat{N}_{A_h}(\bt)c_{\bt,i-1}^{k-1} \right).
\end{split}
\end{equation}
Denoting by $ (c_{\bt,1}^k,c_{\bt,2}^k, \ldots,  c_{\bt,M}^k)$ the approximation of $c_{\bt}^k(t)$ on the defined uniform grid in time, 
we obtain the following matrix form
\[
\left(\begin{matrix}
 c_{\bt,1}^k\\c_{\bt,2}^k\\ \vdots \\ c_{\bt,M}^k
\end{matrix}\right)=\widetilde{\mathcal{M}}_{h,\tau}^{-1}(\bt)\widetilde{\mathcal{N}}_{h,\tau}(\bt)
\left( \begin{matrix}
c_{\bt,1}^{k-1}\\c_{\bt,2}^{k-1}\\ \vdots \\ c_{\bt,M}^{k-1}
\end{matrix} 
\right)
\]
where 
\begin{align*}
\hspace{-1.3cm}&\widetilde{\mathcal{M}}_{h,\tau}(\bt)=\\
\hspace{-1.5cm}&\left(
\begin{smallmatrix}
\hspace{-.25cm}\frac{1}{\tau}\widehat{M}_{B_h}(\bt)+\frac{1}{2}\widehat{M}_{A_h}(\bt)& 0&\cdots &0 \\
\hspace{-.2cm}-\frac{1}{\tau}\widehat{M}_{B_h}(\bt)+\frac{1}{2}\widehat{M}_{A_h}(\bt)&\frac{1}{\tau}\widehat{M}_{B_h}(\bt)+\frac{1}{2}\widehat{M}_{A_h}(\bt)&\cdots&0\\
\hspace{-.25cm}\vdots&\ddots &\ddots&\vdots\\
\hspace{-.25cm}0&\cdots &-\frac{1}{\tau}\widehat{M}_{B_h}(\bt)+\frac{1}{2}\widehat{M}_{A_h}(\bt)&\frac{1}{\tau}\widehat{M}_{B_h}(\bt)+\frac{1}{2}\widehat{M}_{A_h}(\bt)
\end{smallmatrix}
\right).
\end{align*}
 We can obtain $\widetilde{\mathcal{N}}_{h,\tau}(\bt)$  as $\widetilde{\mathcal{M}}_{h,\tau}(\bt)$ by only substituting $\widehat{M}(\bt)$ with $\widehat{N}(\bt)$ in the above matrix.
We can immediately define the smoothing factor of the Gauss-Seidel relaxation operator by considering its symbol $\widetilde{\mathcal{S}}_{h,\tau}(\bt)=\widetilde{\mathcal{M}}_{h,\tau}^{-1}(\bt)\widetilde{\mathcal{N}}_{h,\tau}(\bt)$, that is
\begin{align}
\mu=\sup_{\Theta_h\smallsetminus \Theta_{2h}}\left( \rho\left(\widetilde{\mathcal{S}}_{h,\tau}(\bt)\right) \right).
\end{align} 

Now, we present the analysis of the coarse-grid correction method. As before,  we investigate the effect of the coarse-grid correction, $C_h^{2h}$, on the Fourier modes. The  coarse-grid correction operator is given by: 
\begin{align}
C_h^{2h}=I_h-I_{2h}^h(B_{2h}D_t+\Sigma_t A_{2h})^{-1}I_h^{2h}(B_hD_t+\Sigma_t A_h), 
\end{align} 
where $D_t$ and $\Sigma_t$ are operators corresponding to the Crank-Nicolson approach
and $I_h, I_{2h}^h, I_h^{2h}$ are the identity operator, the transfer operator from coarse to fine grids and vice versa. Also, $(B_{2h}D_t+\Sigma_t A_{2h})$ and $(B_hD_t+\Sigma_t A_h)$ are the coarse- and fine-grid operators, respectively.

Our analysis has to take into account the fact that some of the Fourier modes $\varphi_h(\bt,\cdot)$ on the fine grid coincide on $\mathcal{Q}_{2h}$. For any low frequency $\bt^{00}=(\theta_1 , \theta_2)\in\Theta_{2h}$, we consider the frequencies $\bt^{\alpha}=\bt^{00}-(\alpha_1 \mathrm{sign}(\theta_1),\alpha_2 \mathrm{sign}(\theta_2))\frac{\pi}{h}$, where  $\alpha=\{(\alpha_1,\alpha_2)\vert \alpha_j\in\{0,1\}, j=1,2\}$.
The corresponding four Fourier modes $\varphi_h(\bt^{\alpha},\cdot)$ are called harmonics of each other and they form the space of $2h- $ harmonics for $\bt=\bt^{00}\in\Theta_{2h}$, as follows
\[
\mathcal{F}_{2h}(\bt)=\mathrm{span}\{\varphi_h(\bt^{00},\cdot), \varphi_h(\bt^{11},\cdot),\varphi_h(\bt^{10},\cdot),\varphi_h(\bt^{01},\cdot)\}.
\]
The space of $2h-$harmonics is invariant under the coarse-grid correction operator. Then, the representation of $C_h^{2h}$ on the space $\mathcal{F}_{2h}$ is a $4\times 4$ matrix $\widehat{C}_h^{2h}(\bt).$ To be more precise, we define two vectors: $\boldsymbol{\varphi}_h(\bt,\cdot)=(\varphi_h(\bt^{00},\cdot),$ $\varphi_h(\bt^{11},\cdot),\varphi_h(\bt^{10},\cdot),$ $\varphi_h(\bt^{01},\cdot))$ and ${\bf c}_{\bt}^k(t)=(c_{\bt^{00}}^k(t),c_{\bt^{11}}^k(t),$ $c_{\bt^{10}}^k(t),$ $c_{\bt^{01}}^k(t))$, and then, the error at the  $k$-th iteration  will be $e_h^k({\B},t)\!=\!$ $\sum_{\bt\in\Theta_{2h}}{\bf c}_{\bt}^k(t)$ $\boldsymbol{\varphi}_h(\bt,{\B})^T$. After applying the coarse-grid correction operator on this error, we obtain $\sum_{\bt\in\Theta_{2h}}\widehat{C}_h^{2h}(\bt){\bf c}_{\bt}^k(t)\cdot\boldsymbol{\varphi}_h(\bt,\cdot)$, where $\widehat{C}_h^{2h}(\bt)$ is the following  $4\times 4$ matrix
\begin{align*}
\hspace{-.7cm}\widehat{C}_h^{2h}(\bt)=I_4-\widehat{I}_{2h}^h(\bt)(\widehat{B}_{2h}(\bt)D_t+\Sigma_t\widehat{A}_{2h}(\bt))^{-1}\widehat{I}_h^{2h}(\widehat{B}_h(\bt)D_t+\Sigma_t\widehat{A}_h(\bt)),
\end{align*}
where, $I_4$ is the $4\times 4$ identity matrix, $\widehat{A}_{2h}(\bt)$ and $\widehat{B}_{2h}(\bt)$ are $1\times 1$ symbols of the discrete operators on the coarse grid, and the rest of the involved Fourier symbols are given  by:\\
 $\widehat{A}_h(\bt)=diag\left(\widehat{A}_h(\bt^{00}),\widehat{A}_h(\bt^{11}),\widehat{A}_h(\bt^{10}),\widehat{A}_h(\bt^{01})\right),$\\
$\widehat{B}_h(\bt)=diag\left(\widehat{B}_h(\bt^{00}),\widehat{B}_h(\bt^{11}),\widehat{B}_h(\bt^{10}),\widehat{B}_h(\bt^{01})\right),$\\
$\widehat{I}_{2h}^h(\bt)=\left(\widehat{I}_{2h}^h(\bt^{00}), \widehat{I}_{2h}^h(\bt^{11}),\widehat{I}_{2h}^h(\bt^{10}),\widehat{I}_{2h}^h(\bt^{01}) \right)^T$,\\
$\widehat{I}_{h}^{2h}(\bt)=\left(\widehat{I}_{h}^{2h}(\bt^{00}), \widehat{I}_{h}^{2h}(\bt^{11}),\widehat{I}_{h}^{2h}(\bt^{10}),\widehat{I}_{h}^{2h}(\bt^{01}) \right).$\\

Now, by considering the time discretization of operators $D_t$ and $\Sigma_t$, we obtain the following $M\times M$ matrices:
\[
D_t=\frac{1}{\tau}\left( \begin{matrix}
1&0&\cdots&0\\
-1&1&\cdots&0\\
\vdots&\ddots&\ddots&\vdots\\
0&\cdots&-1&1
\end{matrix}\right),
\quad
\Sigma_t=\frac{1}{2}\left( \begin{matrix}
1&0&\cdots&0\\
1&1&\cdots&0\\
\vdots&\ddots&\ddots&\vdots\\
0&\cdots&1&1
\end{matrix}\right),
\]
Thus,  the error after application of the coarse grid correction is given by  $\widetilde{\mathcal{C}}_{h,\tau}^{2h}(\bt){\bf c}_{\bt}^k(t)\cdot\boldsymbol{\varphi}_h(\bt,\cdot)^T$, where the $4M\times 4M$ matrix  $\widetilde{\mathcal{C}}_{h,\tau}^{2h}(\bt)$ is:
\begin{align*}
\widetilde{\mathcal{C}}_{h,\tau}^{2h}(\bt)=
I_{4M}-\widetilde{\mathcal{I}}_{2h}^h(\bt)\left( \widetilde{\mathcal{K}}_{2h,\tau}(\bt)\right)^{-1}
\widetilde{\mathcal{I}}_{h}^{2h}(\bt)\widetilde{\mathcal{K}}_{h,\tau}(\bt),
\end{align*}
such that $I_{4M}$ is the identity matrix of order $4M$ and $\widetilde{\mathcal{K}}_{h,\tau}(\bt)$ is a  $4M\times 4M$ matrix defined as:
\[
\widetilde{\mathcal{K}}_{h,\tau}(\bt)=\left(
\begin{matrix}
\widetilde{\mathcal{K}}_{h,\tau}(\bt^{00})&0&0&0\\
0&\widetilde{\mathcal{K}}_{h,\tau}(\bt^{11})&0&0\\
0&0&\widetilde{\mathcal{K}}_{h,\tau}(\bt^{10})&0\\
0&0&0&\widetilde{\mathcal{K}}_{h,\tau}(\bt^{01})
\end{matrix}
\right),
\]
with the following blocks in the diagonal:
\begin{align*}
&\hspace{-1cm}\widetilde{\mathcal{K}}_{h,\tau}(\bt^{\alpha})=\\
&\hspace{-1cm}\frac{1}{2\tau}\left(
\begin{smallmatrix}
2\widehat{B}_h(\bt^{\alpha})+\tau\widehat{A}_h(\bt^{\alpha})&&\cdots&0\\
-2\widehat{B}_h(\bt^{\alpha})+\tau\widehat{A}_h(\bt^{\alpha})&2\widehat{B}_h(\bt^{\alpha})+\tau\widehat{A}_h(\bt^{\alpha})&\cdots&0\\
\vdots&\ddots&\ddots&\vdots\\
0&\cdots&-2\widehat{B}_h(\bt^{\alpha})+\tau\widehat{A}_h(\bt^{\alpha})&2\widehat{B}_h(\bt^{\alpha})+\tau\widehat{A}_h(\bt^{\alpha})
\end{smallmatrix}
\right),
\end{align*}
where $\alpha=\{(\alpha_1,\alpha_2)\vert \alpha_j\in\{0,1\}, j=1,2\}$.
In a similar way, we can obtain the Fourier representation of the  prolongation and restriction operators, that are $4M\times M$ and $M\times 4M$ matrices, respectively,
\[\widetilde{\mathcal{I}}_{2h}^h(\bt)=\left( \widehat{I}_{2h}^{h}(\bt^{00})I_M,\, \widehat{I}_{2h}^{h}(\bt^{11})I_M,\,\widehat{I}_{2h}^{h}(\bt^{10})I_M,\,\widehat{I}_{2h}^{h}(\bt^{01})I_M\right)^T,\]
\[\widetilde{\mathcal{I}}_{h}^{2h}(\bt)=\left( \widehat{I}_{h}^{2h}(\bt^{00})I_M,\, \widehat{I}_{h}^{2h}(\bt^{11})I_M,\,\widehat{I}_{h}^{2h}(\bt^{10})I_M,\,\widehat{I}_{h}^{2h}(\bt^{01})I_M\right).\]

Now, we are ready to perform the semi-algebraic two-grid analysis. We do this by combining the presented Fourier smoothing analysis and the Fourier analysis for the coarse-grid correction operator. The two-grid operator is defined as $\mathcal{T}_{h,\tau}^{2h}=\mathcal{S}_{h,\tau}^{\nu_2}\mathcal{C}_{h,\tau}^{2h}\mathcal{S}_{h,\tau}^{\nu_1}$, where, $\mathcal{S}_{h,\tau}$ is the smoothing operator such that the number of pre- and post-smoothing iterations are defined  by $\nu_1$ and $\nu_2$, respectively,
 and $\mathcal{C}_{h,\tau}^{2h}$ is the coarse-grid correction operator.
 
The invariance property of the two-grid method comes from the invariance property of both components of this method. To be more precise,  the coarse-grid correction operator, $\mathcal{C}_{h,\tau}^{2h}$, and the considered Gauss-Seidel smoothing operator, $\mathcal{S}_{h,\tau}$,  both leave the space of $2h-$harmonics $\mathcal{F}_{2h}(\bt^{00})$ invariant for an arbitrary Fourier frequency $\bt^{00}\in \Theta_{2h}$.

 Let us assume again that the error at the $k$th iteration is given as ${\bf c}_{\bt}^k(t)\cdot\boldsymbol{\varphi}_h(\bt,\cdot)^T$. By using the discretization of operators $D_t$ and $\Sigma_t$ we can obtain the relation $\widetilde{\mathcal{T}}_{h,\tau}^{2h}(\bt){\bf c}_{\bt}^k(t)\cdot\boldsymbol{\varphi}_h(\bt,\cdot)^T$ that is the error after the application of the two-grid method, with $\widetilde{\mathcal{T}}_{h,\tau}^{2h}(\bt)$ a  $4M\times 4M$ matrix, given by
\[
\widetilde{\mathcal{T}}_{h,\tau}^{2h}(\bt)=\widetilde{\mathcal{S}}_{h,\tau}^{\nu_2}(\bt)
\left(I_{4M}-\widetilde{\mathcal{I}}_{2h}^h(\bt)\left( \widetilde{\mathcal{K}}_{2h,\tau}(\bt)\right)^{-1}
\widetilde{\mathcal{I}}_{h}^{2h}(\bt)\widetilde{\mathcal{K}}_{h,\tau}(\bt)\right)
\widetilde{\mathcal{S}}_{h,\tau}^{\nu_1}(\bt).
\]
By considering the Gauss-Seidel smoothing operator, the structure of matrix $\widetilde{\mathcal{S}}_{h,\tau}(\bt)$
 is as follows
\[
\widetilde{\mathcal{S}}_{h,\tau}(\bt)=\left(
\begin{matrix}
\widetilde{\mathcal{S}}_{h,\tau}(\bt^{00})&0&0&0\\
0&\widetilde{\mathcal{S}}_{h,\tau}(\bt^{11})&0&0\\
0&0&\widetilde{\mathcal{S}}_{h,\tau}(\bt^{10})&0\\
0&0&0&\widetilde{\mathcal{S}}_{h,\tau}(\bt^{01})
\end{matrix}
\right),
\]
 where the matrix $\widetilde{\mathcal{S}}_{h,\tau}(\bt^{\alpha})$, for $\alpha=\{(\alpha_1,\alpha_2)\vert \alpha_j\in\{0,1\}, j=1,2\}$,  has been previously described in detail.
Finally, we can define the estimation of the two-grid convergence factor by the following expression:
\begin{align}\label{e2.8}
\rho=\sup_{\bt \in \Theta_{2h}}\left( \rho\left(   \widetilde{\mathcal{T}}_{h,\tau}^{2h} (\bt) \right)\right).
\end{align}
\subsection{SAMA results in two dimensions}\label{sec2.2}
In this section we apply the SAMA analysis previously introduced to our model problem. We consider a bilinear finite element discretization on a uniform rectangular mesh to obtain the following discrete problem,
\begin{align*}
B_h\dot{\bu}_h(t)+A_h{\bu}_h(t)=F_h(t),\quad {\bu}_h(0)=g_h,\quad t>0.
\end{align*}
When dealing with finite element discretizations of PDEs, the large sparse stiffness and mass matrices are typically built by the standard assembly procedure \cite{6, karakoc, Dehghan}. In the case of dealing with a structured grid, however, it suffices to represent the discrete operators by means of stencils \cite{stencil}.    The corresponding stencils  obtained for the  mass and stiffness matrices in two dimensions are as follows:
\[
B_h=\frac{h^2}{36}\left[\begin{matrix}
  1 & 4  & 1\\
  4 & 16 & 4\\
   1& 4  & 1 
\end{matrix}\right],\qquad 
A_h=\frac{1}{3}\left[\begin{matrix}
  -1 & -1 & -1\\
  -1 &  8 & -1\\
  -1 & -1 & -1
\end{matrix}\right].
\]
To apply  the presented analysis  we also need to define the components of the multigrid waveform relaxation method. As we stated in the previous sections a Gauss-Seidel waveform relaxation is considered. Regarding the intergrid transfer operators, the  stencil of the restriction operator, $I_{h}^{2h}$, is given by
\[
I_{h}^{2h}=\frac{1}{16}\left[ 
\begin{matrix}
1 & 2 & 1\\
2 & 4 & 2\\
1 & 2 & 1
\end{matrix}
\right].
\]
The prolongation operator ${I}_{2h}^{h}$, is obtained  according to the relation $I_{h}^{2h}=\frac{1}{4}I_{2h}^{h}$, see \cite{2,3} for more details.

\begin{figure}[!h]
\centering
	\centerline{\includegraphics[width=14.5cm]{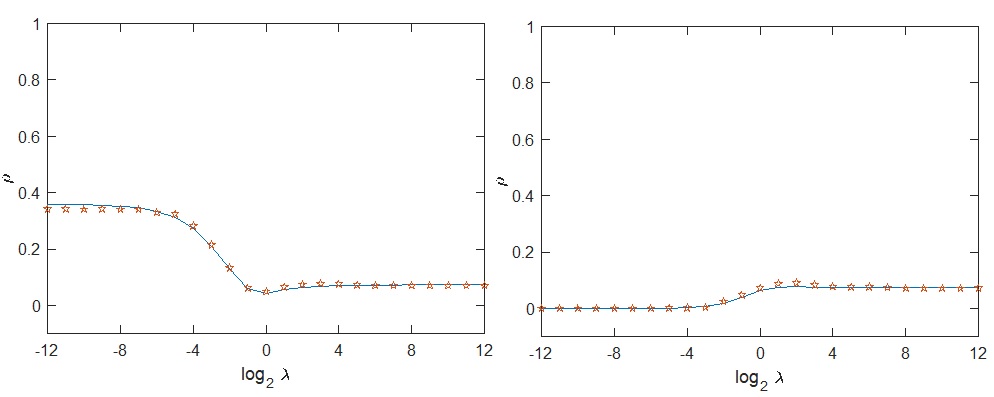}}
\caption{Two-grid convergence factors predicted by the SAMA  together with the experimentally computed convergence factors by W(1,1)-cycle for various values of parameter $\lambda=\tau/h^2$.  The matrix $B_h$ is considered as the mass matrix and the identity matrix in the left and the right pictures, respectively.  }\label{lamb}
\end{figure}

In Figure \ref{lamb}, we show the obtained results by using the semi-algebraic mode analysis { together with the experimentally computed convergence factors by W(1,1)-cycle}. We analyze the convergence of the presented method depending on the parameter $\lambda=\tau/h^2$, which represents the anisotropy in the operator. As we can see in the left picture in  Figure \ref{lamb},  the two-grid convergence factor for various values of parameter $\lambda$  ranging from $2^{-12}$ to $2^{12}$ is bounded by $0.3117$ when matrix $B_h$ is considered as the mass matrix, as stated before. When we consider $B_h$ as the identity matrix, however, the upper bound  is $0.0695$, see the right picture in Figure \ref{lamb}. One can see that the results when choosing $B_h$ as the identity matrix  match with the case in which  the applied discretization method is the finite difference scheme \cite{4}.  In this Figure,  two smoothing steps  are considered and the number of time steps are set to 32 ($M=32$). We should notice that this two-grid  analysis predicts the asymptotic convergence of the so-called W-cycle with two smoothing steps.  Finally, we would like to mention that the asymptotic convergence factor of the proposed multigrid waveform relaxation is very satisfactory concerning the case where $h$ is very small ( $\tau=10^{-5}$ and $\lambda=2^{12}$), where the most iterative procedures lose out.
 \section{SAMA in three dimensions}
 \subsection{Theoretical framework }
We have to do some changes on the analysis to extend the SAMA from two dimensions to three dimensions. As we did for the two dimensional case we need to define  first the infinite grid $\mathcal{Q}_h$,
\[
\mathcal{Q}_h=\{ {\bf x=kh}=(k_1h_1,k_2h_2,k_3h_3), \,\, {\bf k}\in\mathbb{Z}^3\}.
\] 
 In this case, the Fourier modes are $\varphi_h(\bt, {\B})=e^{i\bt\cdot {\B}}=e^{i\theta_1x_1}e^{i\theta_2x_2}e^{i\theta_3x_3}$ where $\bt \in \Theta_h=(-\pi/h,\pi/h]^3$. As in the two dimensional case,  any discrete grid function for a fixed t can be defined as a formal linear combination of Fourier modes . Also, we obtain a new $8-$dimensional Fourier space generated by one low-frequency Fourier mode $\bt=\bt^{000} \in \Theta_{2h}=(-\pi/2h,\pi/2h]^3$ as follows:
\begin{align*}
\mathcal{F}_{2h}(\bt)=&\mathrm{span}\{\varphi_h(\bt^{000},\cdot), \varphi_h(\bt^{111},\cdot),\varphi_h(\bt^{100},\cdot),\\
                      &\varphi_h(\bt^{011},\cdot),\varphi_h(\bt^{010},\cdot),\varphi_h(\bt^{101},\cdot),\varphi_h(\bt^{001},\cdot),\varphi_h(\bt^{110},\cdot)\},
\end{align*}
such that the high-frequencies are given by the following equation
 \[
 \bt^{\alpha}=\bt^{000}-\left(\,\alpha_1 \mathrm{sign}(\theta_1),\alpha_2 \mathrm{sign}(\theta_2),\alpha_3 \mathrm{sign}(\theta_3)\,\right)\,\frac{\pi}{h},
 \]
  where, $\bt^{000}=(\theta_1, \theta_2, \theta_3)$ and $\alpha=\{(\alpha_1,\alpha_2,\alpha_3)\,\vert \, \alpha_j\in\{0,1\}, j=1,2,3\}$.

By the above definitions the resulting Fourier representations of the smoothing, coarse-grid correction and two-grid operators are $8M\times 8M$ matrices that can be computed following the same idea carried out in Section \ref{s2.1} for the two dimensional case. So in fact, SAMA in three dimensions is based on a three dimensional spatial LFA combined with an exact analysis in the time variable. 
\subsection{SAMA results in three dimensions}
Here we consider the discretization of our model problem by using trilinear finite elements. Similarly to the two dimensional case,  first we define the stencil forms of the mass and stiffness matrices. By some computations we can obtain the stencils for  $B_h$ and $A_h$  in three dimensions as follows:
\[
B_h=\frac{h^3}{216}\left[
\left[\begin{matrix}
  1 & 4  & 1\\
  4 & 16 & 4\\
   1& 4  & 1 
\end{matrix}\right]     
\left[\begin{matrix}
  4  & 16 & 4\\
  16 & 64 & 16\\
   4 & 16 & 4 
\end{matrix}\right]  
\left[\begin{matrix}
  1 & 4  & 1\\
  4 & 16 & 4\\
  1 & 4  & 1 
\end{matrix}\right]  
\right],
\]  \[
A_h= \frac{h}{12}\left[
\left[\begin{matrix}
  -1 & -2 & -1\\
  -2 &  0 & -2\\
   -1& -2 & -1 
\end{matrix}\right]     
\left[\begin{matrix}
  -2 & 0  & -2\\
   0 & 32 & 0\\
  -2 & 0  & -2 
\end{matrix}\right]  
\left[\begin{matrix}
 -1 & -2 & -1\\
  -2 &  0 & -2\\
   -1& -2 & -1 
\end{matrix}\right]  
\right].
\]

\begin{figure}[!h]
\centering
	\centerline{\includegraphics[width=14.5cm]{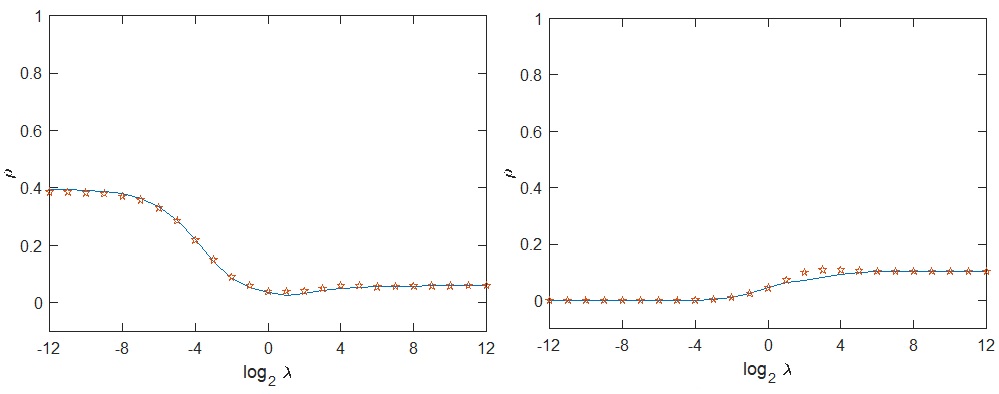}}
\caption{Two-grid convergence factors  predicted by the  SAMA { together with the experimentally computed convergence factors by W(1,1)-cycle} for various values of parameter $\lambda=\tau/h^2$.  Matrix $B_h$ is considered as the mass matrix and the identity matrix  in the left and the right pictures, respectively.  }\label{lamb2}
\end{figure}

Also, we need to fix the multigrid components. As before, we use a Gauss-Seidel waveform relaxation and the inter-grid transfer operators are given as follows. The stencil of the restriction operator, $I_{h}^{2h}$, is,
\[
I_{h}^{2h}=\frac{1}{64}\left[
\left[
\begin{matrix}
1 & 2 & 1\\
2 & 4 & 2\\
1 & 2 & 1
\end{matrix}\right]
\left[
\begin{matrix}
2 & 4 & 2\\
4 & 8 & 4\\
2 & 4 & 2
\end{matrix}\right]
\left[
\begin{matrix}
1 & 2 & 1\\
2 & 4 & 2\\
1 & 2 & 1
\end{matrix}\right]
\right],
\]
 and the prolongation operator ${I}_{2h}^{h}$ will satisfy  $I_{h}^{2h}=\frac{1}{8}I_{2h}^{h}$.
In Figure \ref{lamb2}, we show the obtained results by two-grid SAMA analysis in three dimensions  together with the experimentally computed convergence factors by W(1,1)-cycle. As in the two dimensional case, we analyze the convergence of the presented method accordingly to  parameter $\lambda=\tau/h^2$ (which denotes the anisotropy in the operator). We consider again two different cases: when $B_h$ is the mass matrix (left picture) and when $B_h$ is assumed to be the identity matrix (right picture). Here we obtain very similar results to the two dimensional case. To be more precise, when matrix $B_h$ is the identity matrix, the results   match very good with the case in which a finite difference discretization method is used. As we expected, the upper bound for the predicted convergence factor slightly increases  for both cases being equal to $0.4040$ and $0.0881$ when $B_h$ is the mass matrix and identity matrix, respectively.
  The results are presented for two smoothing steps  and  $M=32$. {Moreover as expected, the asymptotic convergence factor of the proposed multigrid waveform relaxation is very satisfactory concerning the case where $h$ tends to be very small.} 
\section{Conclusions}
 In this work, we presented the two-grid SAMA analysis for predicting the convergence factor of a multigrid waveform relaxation method for finite element discretizations in two and three spatial dimensions.  We considered the Crank-Nicolson discretization in time and bilinear and trilinear finite element methods for the spatial discretization in two and three dimensions, respectively. The presented results when matrix $B_h$  is the identity matrix, are comparable with the case in which a finite difference discretization method is used. The proposed SAMA analysis allows us to systematically study the behavior of the multigrid waveform relaxation method for finite element discretizations.

\end{document}